\documentclass[a4paper,11pt]{article}
\usepackage{graphicx} 
\usepackage{amsmath}
\usepackage{amssymb}
\usepackage{amsfonts}
\usepackage{amsthm}
\usepackage{parskip}
\usepackage{comment}
\usepackage{xcolor}
\usepackage{enumitem}
\usepackage[hidelinks]{hyperref}
\usepackage[top=22mm, bottom=23mm, left=22mm, right=22mm]{geometry}

\newtheorem{thm}{Theorem}[section]
\newtheorem{prop}[thm]{Proposition}
\newtheorem{lem}[thm]{Lemma}
\newtheorem{cor}[thm]{Corollary}
\newtheorem{conj}[thm]{Conjecture}
\newtheorem{ques}[thm]{Question}

\theoremstyle{definition}

\newtheorem{remark}[thm]{\bf Remark}
\newtheorem{obs}[thm]{\bf Observation}
\newtheorem{df}[thm]{\bf Definition}

\title{A variant of the Erdős-Gyárfás problem for $K_8$} 
\author{Fredy Yip\thanks{Trinity College, University of Cambridge, United Kingdom. Email: \textbf{fy276@cam.ac.uk}.}}
\date{}
\newcommand{\floor}[1]{\left\lfloor #1 \right\rfloor}

\begin{document}

\maketitle

\begin{abstract}
    Recently, Alon initiated the study of graph codes and their linear variants in analogy to the study of error correcting codes in theoretical computer science. Alon related the maximum density of a linear graph code which avoids images of a small graph $H$ to the following variant of the Erdős-Gyárfás problem on edge-colourings of $K_n$. A copy of $H$ in an edge-colouring of $K_n$ is \emph{even-chromatic} if each colour occupies an even number of edges in the copy. We seek an edge-colouring of $K_n$ using $n^{o(1)}$ colours such that there are no even-chromatic copies of $H$. Such an edge-colouring is conjectured to exist for all cliques $K_t$ with an even number of edges. To date, edge-colourings satisfying this property have been constructed for $K_4$ and $K_5$. 

    We construct an edge-colouring using $n^{o(1)}$ colours which avoids even-chromatic copies of $K_8$. This was the smallest open case of the above conjecture, as $K_6, K_7$ each has an odd number of edges. We also study a stronger condition on edge-colourings, where for each copy of $H$, there is a colour occupying exactly one edge in the copy. We conjecture that an edge-colouring using $n^{o(1)}$ colours and satisfying this stronger requirement exists for all cliques $K_t$ regardless of the parity of the number of its edges. We construct edge-colourings satisfying this stronger property for $K_4$ and $K_5$. These constructions also improve upon the number of colours needed for the original problem of avoiding even-chromatic copies of $K_4$ and $K_5$.
\end{abstract}

\section{Introduction}

In theoretical computer science, an error correcting code is a method of encoding information redundantly, so as to allow its recovery in the presence of limited transmission errors. Mathematically, an error correcting code using $n$-bits may be represented as a set $\mathcal{C}\subseteq\mathbb{F}_2^{[n]}$ of binary strings of length $n$, where the Hamming distance between distinct elements of $\mathcal{C}$ is suitably bounded from below. This is equivalent to the absence of sparse non-zero binary strings in $\mathcal{C} + \mathcal{C} = \{s_1 + s_2|s_1, s_2\in\mathcal{C}\}$. 

In \cite{alon2024graph}, Alon introduced an analogous notion of ``graph codes". Here, in place of $n$-bit binary strings, we consider graphs on vertex set $[n]$. We may represent a graph $G$ on vertex set $[n]$ by its edge set $E(G)\subseteq \binom{[n]}{2}$ as an element of $\mathbb{F}_2^{\binom{[n]}{2}}$, where edges correspond to $1$-entries and non-edges correspond to $0$-entries. For two graphs $G_1, G_2$ on vertex set $[n]$, the addition operation on $\mathbb{F}_2^{\binom{[n]}{2}}$ gives a graph $G_1 + G_2$ on vertex set $[n]$ whose edge set is given by the \emph{symmetric difference} of the edge sets $E(G_1)$ and $E(G_2)$. 

The analogous notion of sparse binary strings in the graph setting is given by the \emph{image} of certain small graphs. Here an \emph{image} of a graph $H$ (not necessarily on vertex set $[n]$) is a graph $G$ on vertex set $[n]$ consisting solely of an induced subgraph isomorphic to $H$, with no other edges. In other words, for a subset $S\subseteq [n]$ of vertices, $G[S]$ is isomorphic to $H$ and all vertices in $[n]\backslash S$ are isolated. 

Formally, we have the following definition of a graph code with respect to a family $\mathcal{H}$ of forbidden graphs. 

\begin{df}
    For a collection $\mathcal{H}$ of graphs, an $\mathcal{H}$-(graph)-code on vertex set $[n]$ is a collection $\mathcal{C}\subseteq \mathbb{F}_2^{\binom{[n]}{2}}$ of graphs on vertex set $[n]$ for which $\mathcal{C} + \mathcal{C} = \{G_1 + G_2|G_1, G_2\in\mathcal{C}\}$ does not contain images of any $H\in\mathcal{H}$. 
\end{df}

In analogy to the efficiency of error correcting codes, Alon et al.~\cite{alon2024graph, alon2023structured} investigated the maximum density of an $\mathcal{H}$-code on vertex set $[n]$, defined as
\begin{equation*}
    d_{\mathcal{H}}(n) = \sup_{\substack{\mathcal{H}\text{-code }\mathcal{C}\text{ on}\\\text{vertex set } [n]}} \frac{|\mathcal{C}|}{2^{\binom{n}{2}}}. 
\end{equation*}

For any family $\mathcal{H}$ of forbidden graphs, $d_{\mathcal{H}}(n)$ is a decreasing function of $n$. A more general problem, where $\mathcal{C} + \mathcal{C}$ must avoid a family of graphs \emph{on vertex set $[n]$} (\textit{i.e.}~without considering images) has been studied in \cite{alon2023structured}. In particular, \cite{alon2023structured} considered the family of connected graphs, the family of 2-connected graphs, the family of Hamiltonian graphs, the family of graphs containing a spanning star and the complements of each of these families. The case of graph-codes was also considered in the same work under the framework of families on vertex set $[n]$ defined by local conditions. For the family $\mathcal{K}$ of all cliques, a conjecture of Gowers \cite{gowers} implies that $d_{\mathcal{K}}(n) = o(1)$, which remains open. 

One may also consider the linear setting, where the graph code $\mathcal{C}\leq \mathbb{F}_2^{\binom{[n]}{2}}$ is a vector subspace. In this setting we have $\mathcal{C} + \mathcal{C} = \mathcal{C}$, resulting in the following definition of linear graph codes. 

\begin{df}
    For a collection $\mathcal{H}$ of graphs, a linear $\mathcal{H}$-(graph)-code on vertex set $[n]$ is a vector subspace $\mathcal{C}\leq \mathbb{F}_2^{\binom{[n]}{2}}$ of graphs on vertex set $[n]$ for which $\mathcal{C}$ does not contain images of any $H\in\mathcal{H}$. 
\end{df}

We may similarly consider the maximum density of a linear graph code, defined as

\begin{equation*}
    d^\text{lin}_{\mathcal{H}}(n) = \sup_{\substack{\text{linear }\mathcal{H}\text{-code }\mathcal{C}\\\text{on vertex set } [n]}} \frac{|\mathcal{C}|}{2^{\binom{n}{2}}}. 
\end{equation*}

In general, we have $d^\text{lin}_{\mathcal{H}}(n)\leq d_{\mathcal{H}}(n)$. Like $d_{\mathcal{H}}(n)$, $d^{\text{lin}}_{\mathcal{H}}(n)$ is a decreasing function of $n$ for any family $\mathcal{H}$ of forbidden graphs. 

The special case where $\mathcal{H} = \{H\}$ consists of a single graph is of key interest. Here we use the notation $d_H(n) = d_{\{H\}}(n)$ and $d^{\text{lin}}_H(n) = d^{\text{lin}}_{\{H\}}(n)$. For a graph $H$ with an odd number of edges, we have $d_H(n) = d^{\text{lin}}_H(n) = 1/2$ for all sufficiently large~$n$, attained by the linear graph code consisting of all graphs on $[n]$ with an even number of edges. Therefore, we restrict our attention to the case where $H$ has an even number of edges. 

Alon \cite{alon2024graph} determined the order of magnitude of $d_H(n)$ when $H$ is a star or a matching. Conlon, Lee and Versteegen \cite{conlon2024around} proved that $d_H(n)=o(1)$ when $H$ is a positive graph. Interestingly, as noted in \cite{versteegen2023upper}, all known lower bounds to $d_H(n)$ are given via linear graph codes. 

Alon \cite{alon2024graph} related the asymptotic behaviour of $d^{\text{lin}}_{H}(n)$ to the following Ramsey-theoretic problem. 

\begin{df}
    Given an edge-colouring of $K_n$, a copy of $H$ in $K_n$ is \emph{even-chromatic} if it has an even number of edges of each colour. An edge-colouring of $K_n$ is \emph{$H$-odd} if it does not admit even-chromatic copies of $H$. Let $r_H(n)$ be the least number of colours needed for an $H$-odd edge-colouring of $K_n$. 
\end{df}

As with $d_H(n)$ and $d_H^\text{lin}(n)$, the behaviour of $r_H(n)$ is only of interest when $H$ has an even number of edges, as $r_H(n) = 1$ otherwise. It can be shown (see \cite{versteegen2023upper} for a proof) that for any graph $H$ there exist constants $c$ and $C$ (dependent only on $H$) such that
\begin{equation} \label{lin graph code and even chrom relation}
    c\cdot\frac{1}{r_H(n)^{C}}\leq d^\text{lin}_{H}(n)\leq \frac{1}{r_H(n)}. 
\end{equation}

For any graph $H$ with an even number of edges, Versteegen \cite{versteegen2023upper} gave a lower bound of $r_H(n) = \Omega(\log{n})$ with an absolute implied constant. 

We always have $r_H(n)\leq \binom{n}{2}$ by using distinct colours for each edge in $K_n$. Therefore, for any graph~$H$, (\ref{lin graph code and even chrom relation}) yields the polynomial lower bound
\begin{equation*}
    d_H(n)\geq d^{\text{lin}}_H(n) = \Omega\left(n^{-C_H}\right), 
\end{equation*}
for some constant $C_H > 0$. It is of key interest to decide for which graphs we have $d_H(n)\geq d^{\text{lin}}_H(n) = n^{-o(1)}$. 

\begin{ques} \label{lin code subpoly}
    For which graphs $H$, with an even number of edges, do we have $d^{\text{lin}}_H(n) = n^{-o(1)}$? 
\end{ques}

By virtue of (\ref{lin graph code and even chrom relation}), Question \ref{lin code subpoly} is equivalent to the following sub-polynomial upper bound for $r_H(n)$. 

\begin{ques} \label{ramsey subpoly}
    For which graphs $H$, with an even number of edges, do we have $r_H(n) = n^{o(1)}$? 
\end{ques}

Versteegen \cite{versteegen2023upper} defined a class of \emph{even-decomposable} graphs for which he gave a negative answer to Question \ref{lin code subpoly}/\ref{ramsey subpoly}.  

\begin{df}
    A graph $H$ with an even number of edges is \emph{even-decomposable} if there is a sequence $V(H) = V_0\supset V_1\supset\cdots \supset V_k = \emptyset$ such that for each $0\leq i\leq k-1$, $H[V_i]$ has an even number of edges and $V_{i}\backslash V_{i+1}$ is an independent set in~$H$.
\end{df}

\begin{thm}[Versteegen \cite{versteegen2023upper}] \label{versteegen thm}
    If $H$ is even-decomposable, then $r_H(n) = n^{\Omega(1)}$. 
\end{thm}

Janzer and Yip \cite{janzerandyip} showed that the probability that a uniform random graph on $t$ vertices with an even number of edges is even-decomposable is $1 - e^{-\Theta(t^2)}$. Versteegen conjectured that being even-decomposable is the only obstruction to a positive answer to Question \ref{lin code subpoly}/\ref{ramsey subpoly}. 

\begin{conj}[Versteegen \cite{versteegen2023upper}] \label{versteegen conj}
    We have $r_H(n) = n^{\Omega(1)}$ if and only if $H$ is even-decomposable. 
\end{conj}

We shall focus on the case of cliques $H = K_t$. The requirement that $H$ has an even number of edges is equivalent to $t \equiv 0, 1\mod{4}$. It can be seen that for all $t\geq 4$ satisfying this modular condition, $K_t$ is not even-decomposable. Therefore, Conjecture \ref{versteegen conj} implies that $r_{K_t}(n) = n^{o(1)}$ for any $t\geq 4$ with $t \equiv 0, 1\mod{4}$. This special case of Conjecture \ref{versteegen conj} has been conjectured earlier in \cite{ge2023new} by Ge, Xu and Zhang, and remains open. 

\begin{conj}[Ge, Xu and Zhang \cite{ge2023new}] \label{clique n^o(1) conj}
    For any positive integer $t\geq 4$ with $t \equiv 0, 1\mod{4}$, we have
    \begin{equation*}
        r_{K_t}(n) = n^{o(1)}. 
    \end{equation*}
\end{conj}

Affirmative answers were given to Conjecture \ref{clique n^o(1) conj} in the cases of $t = 4$ (Cameron and Heath \cite{cameron2023new}) and $t = 5$ (Ge, Xu and Zhang \cite{ge2023new} and Bennett, Heath and Zerbib \cite{heath2023edge}, independently). 

\begin{thm}[Cameron and Heath \cite{cameron2023new}] \label{thm:CH}
    We have 
    \begin{equation*}
        r_{K_4}(n) = e^{O((\log^{1/2}n)\log\log n)} = n^{o(1)}. 
    \end{equation*}
\end{thm}

\begin{thm}[Bennett, Heath and Zerbib \cite{heath2023edge}; Ge, Xu and Zhang \cite{ge2023new}] \label{thm:GXZ}
    We have 
    \begin{equation*}
        r_{K_5}(n) = e^{O((\log^{1/2}n)\log\log n)} = n^{o(1)}. 
    \end{equation*}
\end{thm}

We remark that Ge, Xu and Zhang obtained the weaker upper bound of $r_{K_5}(n)=e^{O((\log^{2/3}n)\log\log n)}$ and Versteegen \cite{versteegen2023upper} showed that the $\log\log n$ factor may be removed from the exponent of the upper bound for $r_{K_4}(n)$. 

We give an affirmative answer to Conjecture \ref{clique n^o(1) conj} for the smallest open case $t = 8$. 

\begin{thm} \label{main thm}
    We have 
    \begin{equation*}
        r_{K_8}(n) = e^{O(\log^{2/3} n)} = n^{o(1)}. 
    \end{equation*}
\end{thm}

Theorem \ref{main thm} gives, via (\ref{lin graph code and even chrom relation}), the following lower bound on the maximum density of a (linear) $K_8$-code. 

\begin{cor}
    We have 
    \begin{equation*}
        d_{K_8}(n) \geq d^{\text{lin}}_{K_8}(n) = e^{-O(\log^{2/3} n)} = n^{-o(1)}. 
    \end{equation*}
\end{cor}

We shall construct our colouring inductively, where given a $K_8$-odd colouring of $K_n$, we produce a $K_8$-odd colouring of $K_{nm}$ for an appropriate choice of $m = m(n)$. 

We shall first apply our method to the simpler problem of constructing $K_4$-odd and $K_5$-odd colourings with $n^{o(1)}$ colours. Importantly, our colourings will have the stronger property that there is an edge of unique colour in each copy of $K_4$ or $K_5$. We introduce the following strengthening of $H$-oddness to incorporate this stronger property. 

\begin{df}
    Given an edge-colouring of $K_n$, a copy of $H$ in $K_n$ is \emph{unique-chromatic} if there is a colour occupying exactly one edge of the copy. An edge-colouring of $K_n$ is \emph{$H$-unique} if all copies of $H$ in $K_n$ are unique-chromatic. Let $u_H(n)$ be the least number of colours needed for an $H$-unique edge-colouring of $K_n$. 
\end{df}

As any $H$-unique edge-colouring is $H$-odd, we have $r_H(n)\leq u_H(n)$. The author is informed by Conlon that the notions of unique-chromaticity and $H$-uniqueness (for cliques) have been introduced by Radoicic and studied by Axenovich and Conlon in their unpublished works. In particular, Conjecture \ref{clique conj} below had been made by Conlon, to which he gave affirmative answers for $t\leq 7$.

\begin{prop} \label{k4 prop}
    We have
    \begin{equation*}
        r_{K_4}(n)\leq u_{K_4}(n) = e^{O(\log^{1/2} n)} = n^{o(1)}. 
    \end{equation*}
\end{prop}

\begin{prop} \label{k5 prop}
    We have
    \begin{equation*}
        r_{K_5}(n)\leq u_{K_5}(n) = e^{O(\log^{1/2} n)} = n^{o(1)}. 
    \end{equation*}
\end{prop}

Note that Propositions \ref{k4 prop} and \ref{k5 prop} provide alternative proofs of Theorems \ref{thm:CH} and \ref{thm:GXZ}, respectively. Furthermore, Proposition \ref{k5 prop} improves the best known upper bound on $r_{K_5}(n)$ by a factor of $\log\log n$ in the exponent. 

Unlike $H$-oddness, the notion of $H$-uniqueness is meaningful regardless of the parity of the number of edges of $H$. As the addition/removal of isolated vertices from $H$ does not change $H$-oddness/uniqueness, we may assume that $H$ has no isolated vertices. We pose the analogous question to Question \ref{ramsey subpoly} for $H$-uniqueness in place of $H$-oddness. 

\begin{ques} \label{unique subpoly}
    For which graphs $H$, without isolated vertices, do we have $u_H(n) = n^{o(1)}$? 
\end{ques}

Versteegen's \cite{versteegen2023upper} method of proving Theorem \ref{versteegen thm} may be adapted to give a negative answer to Question~\ref{unique subpoly} whenever $H$ is not a clique. (A proof is given in Appendix \ref{non-clique implies non-unique appendix}.)

\begin{prop} \label{non-clique implies non-unique}
    For a non-complete graph $H$ without isolated vertices, we have
    \begin{equation*}
        u_H(n) = \Omega\left(n^{\frac{1}{|V(H)| - 1}}\right) = n^{\Omega(1)}. 
    \end{equation*}
\end{prop}

We conjecture that the converse is true. 

\begin{conj} \label{clique conj}
    For any positive integer $t\geq 2$, we have $u_{K_t}(n) = n^{o(1)}$. 
\end{conj}

As $r_H(n)\leq u_H(n)$, Conjecture \ref{clique conj} implies Conjecture \ref{clique n^o(1) conj}. The Erdős-Gyárfás problem \cite{erdos75, erdos81, erdosgyarfas97} on generalised Ramsey numbers asks for the least number $f_{s, t}(n)$ of colours needed to colour the edges of $K_n$ such that each copy of $K_t$ receives at least $s$ colours. It has been proven in \cite{erdosgyarfasproblem} that $f_{t - 1, t}(n) = n^{o(1)}$. Conjecture \ref{clique conj}, if true, offers another proof of this result. Indeed, taking $c_p$ to be a $K_p$-unique edge-colouring of $K_n$ using $n^{o(1)}$ colours, the product colouring of $c_2, \cdots, c_t$ would suffice. In fact, despite the difference in presentation, the colourings we construct share the multi-levelled structure of the colouring in \cite{erdosgyarfasproblem}. 

\textbf{Organisation of the paper.} In Section \ref{general}, we discuss our general method of inductive colouring constructions via \emph{amalgamations} and bounds on the number of colours used. We make general observations about $H$-oddness and $H$-uniqueness in Section \ref{elem}. In Section \ref{k4}, we apply this method to produce a $K_4$-unique edge-colouring of $K_n$ using $n^{o(1)}$ colours, proving Proposition \ref{k4 prop}. In Section \ref{k5}, we do the same for $K_5$, proving Proposition \ref{k5 prop}. In Section \ref{k8}, we prove Theorem \ref{main thm}, constructing a $K_8$-odd edge-colouring of $K_n$ using $n^{o(1)}$ colours. Appendix \ref{colour counting appendix} proves the technical results giving the $n^{o(1)}$ upper bound on the number of colours used. In Appendix \ref{non-clique implies non-unique appendix} we prove Proposition \ref{non-clique implies non-unique} giving the polynomial lower bound $n^{\Omega(1)}$ for $u_H(n)$ in the case of a non-complete graph $H$ without isolated vertices. 

\textbf{Notation.} For a colouring $c$, we denote by $\mathfrak{p}(c)$ the set of colours used by $c$ and we denote by $\mathfrak{c}(c) = |\mathfrak{p}(c)|$ the number of colours used by $c$. For a graph $G$, we denote by $V(G)$ and $E(G)$ the vertex set and edge set of $G$, respectively. For a subset $S\subseteq V(G)$ of vertices, let $E(S) = E(G[S])$ be the set of edges of $G$ between vertices in $S$. For disjoint subsets $S_1, S_2\subseteq V(G)$ of vertices, we denote by $E(S_1, S_2)$ the set of edges of $G$ between $S_1$ and $S_2$. For an edge-colouring $c: E(K_n)\rightarrow \mathfrak{p}(c)$ of $K_n$, we denote by $c(e) = c(u, v)$ the colour assigned to the edge $e = uv$. We identify a copy of $K_t$ in $K_n$ with its vertex set $S\subseteq V(K_n)$. The elements of a Cartesian product are canonically denoted by appropriate ordered tuples. Given a set $S$, we denote by $S^+ = S\sqcup\{*\}$ the extension of $S$ by adding a new element $*$, which we consider to be distinct from the elements of $S$. We denote by $\log$ the natural logarithm with base $e$. 

\section{Inductive colouring construction} \label{general}

For a property $\mathcal{P}$ (\emph{e.g.}~$K_4$-uniqueness) of edge-colourings of $K_n$, we call an edge-colouring of $K_n$ satisfying $\mathcal{P}$ \emph{a $\mathcal{P}$-colouring of $K_n$}. 

\begin{df}
    For a property $\mathcal{P}$ of edge-colourings of $K_n$, let $\mathfrak{r}_{\mathcal{P}}(n)$ be the least number of colours needed for a $\mathcal{P}$-colouring of $K_n$. 
\end{df}

As examples, we have $r_H(n) = \mathfrak{r}_{H\text{-odd}}(n)$ and $u_H(n) = \mathfrak{r}_{H\text{-unique}}(n)$. 

Given $\mathcal{P}$, we shall construct a $\mathcal{P}$-colouring of $K_n$ by induction on $n$. More precisely, given a $\mathcal{P}$-colouring of $K_n$, we aim to produce a $\mathcal{P}$-colouring of $K_{nm}$, for an appropriate choice of $m = m(n)$, with the help of an auxiliary edge-colouring of $K_m$ satisfying an often ``weaker" property $\mathcal{Q}$. 

To do so, we make use of the following method of \emph{amalgamating} an edge-colouring $c$ of $K_n$ with an auxiliary edge-colouring $d$ of $K_m$ to produce an edge-colouring of $K_{nm}$, which we denote by $c\otimes d$. Here we take the vertex sets of $K_n, K_m, K_{nm}$ to be $[n], [m]$ and the Cartesian product $[n]\times [m]$ respectively. We shall visualise the vertex set $[n]\times [m]$ in a grid with horizontal axis $[n]$ and vertical axis $[m]$. 

\begin{df}
    For an edge-colouring $c$ of $K_n$ (on vertex set $[n]$) and an edge-colouring $d$ of $K_m$ (on vertex set $[m]$), the \emph{amalgamation of $c$ and $d$} is an edge-colouring $c\otimes d$ of $K_{nm}$ (on vertex set $[n]\times [m]$). We shall take the co-domain of $c\otimes d$ to be the $4$-fold Cartesian product $\mathfrak{p}(c)^+\times\mathfrak{p}(d)^+\times\{+, -, 0, \infty\}\times\binom{[m]}{2}^+$. For distinct vertices $(x_1, y_1), (x_2, y_2)$ in $[n]\times [m]$, where $x_1\leq x_2$, we take
    \begin{equation*}
        c\otimes d((x_1, y_1), (x_2, y_2)) = \begin{cases}
            (c(x_1, x_2), d(y_1, y_2), +, *) &\text{ if } x_1 < x_2 \text{ and } y_1 < y_2, \\
            (c(x_1, x_2), d(y_1, y_2), -, *) &\text{ if } x_1 < x_2 \text{ and } y_1 > y_2, \\
            (*, d(y_1, y_2), \infty, \{y_1, y_2\}) &\text{ if } x_1 = x_2, \\
            (c(x_1, x_2), *, 0, *) &\text{ if } y_1 = y_2. \\
        \end{cases}
    \end{equation*}
    We denote the $4$ components of $c\otimes d$ as $(c\otimes d)_1, (c\otimes d)_2, (c\otimes d)_3, (c\otimes d)_4$. 
\end{df}

The first component $(c\otimes d)_1$ applies the edge-colouring $c$ to the edge $x_1x_2$ obtained by projecting the edge $(x_1, y_1)(x_2, y_2)$ onto the horizontal axis $[n]$. This is not well-defined when $x_1 = x_2$, and in this case $(c\otimes d)_1$ assigns the new colour $*$ (distinct from the colours $\mathfrak{p}(c)$ used by $c$) to the edge $(x_1, y_1)(x_2, y_2)$. The second component $(c\otimes d)_2$ applies the edge-colouring $d$ to the edge $y_1y_2$ obtained by projecting the edge $(x_1, y_1)(x_2, y_2)$ onto the vertical axis $[m]$. A similar treatment as above is given for the case of $y_1 = y_2$. The third component $(c\otimes d)_3$ records the sign of the gradient of the edge $(x_1, y_1)(x_2, y_2)$ viewed as a line in the grid. The fourth component $(c\otimes d)_4$ assigns each `vertical' edge $(x, y_1)(x, y_2)$ a colour labelled by the unordered pair $\{y_1, y_2\}$. 

$c\otimes d$ does not use all colours in its co-domain $\mathfrak{p}(c)^+\times\mathfrak{p}(d)^+\times\{+, -, 0, \infty\}\times\binom{[m]}{2}^+$. 

\begin{lem} \label{amalgamation colour count}
    For any edge-colourings $c$ and $d$, we have
    \begin{equation*}
        \mathfrak{c}(c\otimes d) = \left(2\mathfrak{c}(d) + 1\right)\mathfrak{c}(c) + \binom{m}{2}. 
    \end{equation*}
\end{lem}

\begin{proof}
    We count the colour usage for each of the possible colours $+, -, 0, \infty$ of the third component. 
    \begin{itemize}
        \item The third component is `$+$'. We have $\mathfrak{c}(c)$ choices for the first component and $\mathfrak{c}(d)$ choices for the second component with the fourth component fixed. Hence a total of $\mathfrak{c}(c)\mathfrak{c}(d)$ colours are needed in this case. 
        \item The third component is `$-$'. Similarly, we also need a total of $\mathfrak{c}(c)\mathfrak{c}(d)$ colours in this case. 
        \item The third component is `$\infty$'. We use a colour for each unordered pair $\{y_1, y_2\}$ of distinct elements in $[m]$, thus a total of $\binom{m}{2}$ colours are needed in this case. 
        \item The third component is `$0$'. We have $\mathfrak{c}(c)$ choices for the first component, with the second and the fourth components fixed. Hence a total of $\mathfrak{c}(c)$ colours are needed in this case. 
    \end{itemize}
    Taking the sum of the colour usage in each case gives the desired count of the total number of colours $\mathfrak{c}(c\otimes d)$ used by $c\otimes d$. 
\end{proof}

Given a property $\mathcal{P}$, to upper bound $\mathfrak{r}_{\mathcal{P}}(n)$, our strategy will be to identify a ``weaker" property $\mathcal{Q}$ such that the amalgamation of a $\mathcal{P}$-colouring of $K_n$ and a $\mathcal{Q}$-colouring of $K_m$ always gives a $\mathcal{P}$-colouring of $K_{nm}$. Once we establish this, the formula for the number of colours used by an amalgamation in Lemma \ref{amalgamation colour count} gives
\begin{equation} \label{key amalgamation bound}
    \mathfrak{r}_\mathcal{P}(nm) \leq \left(2\mathfrak{r}_\mathcal{Q}(m) + 1\right)\mathfrak{r}_\mathcal{P}(n) + \binom{m}{2}, 
\end{equation}
for any positive integers $m, n$. We shall choose $m = m(n)$ so that the first term on the right-hand side of (\ref{key amalgamation bound}) dominates. A sub-polynomial upper bound $\mathfrak{r}_\mathcal{Q}(n) = n^{o(1)}$ may then by converted to a sub-polynomial upper bound $\mathfrak{r}_\mathcal{P}(n) = n^{o(1)}$, via repeated amalgamations. The following lemma formulates this result precisely. 

\begin{lem} \label{subpoly counting}
    If the amalgamation of a $\mathcal{P}$-colouring and a $\mathcal{Q}$-colouring is always a $\mathcal{P}$-colouring, $\mathfrak{r}_\mathcal{Q}(n) = n^{o(1)}$ and $\mathfrak{r}_\mathcal{P}(n) < \infty$ is an increasing function of $n$, then we have
    \begin{equation*}
        \mathfrak{r}_\mathcal{P}(n) = n^{o(1)}. 
    \end{equation*}
\end{lem}

The condition that $\mathfrak{r}_\mathcal{P}(n)$ is an increasing function of $n$ is needed as the recursion (\ref{key amalgamation bound}) does not directly bound $\mathfrak{r}_\mathcal{P}(n)$ for every $n$. Similarly, the requirement that $\mathfrak{r}_\mathcal{P}(n) < \infty$ ensures that a $\mathcal{P}$-colouring exists. For any (appropriate) graph $H$, both properties are satisfied for both $H$-oddness and $H$-uniqueness. 

In the case that we have a tighter quasi-polynomial bound $\mathfrak{r}_\mathcal{Q}(n) = e^{O(\log^q n)}$, where $q < 1$ is a constant, we may obtain a tighter quasi-polynomial bound $\mathfrak{r}_\mathcal{P}(n) = e^{O(\log^p n)}$, where $p = \frac{1}{2 - q} < 1$. 

\begin{lem} \label{quasipoly counting}
    If the amalgamation of a $\mathcal{P}$-colouring and a $\mathcal{Q}$-colouring is always a $\mathcal{P}$-colouring, $\mathfrak{r}_\mathcal{Q}(n) = e^{O(\log^q n)}$ for a constant $q\in [0, 1)$ and $\mathfrak{r}_\mathcal{P}(n) < \infty$ is an increasing function of $n$, then we have
    \begin{equation*}
        \mathfrak{r}_\mathcal{P}(n) = e^{O(\log^p n)}, 
    \end{equation*}
    where $p = \frac{1}{2 - q} < 1$. 
\end{lem}

We now prove Lemma \ref{quasipoly counting}. The proof of the unused Lemma \ref{subpoly counting} is deferred to Appendix \ref{colour counting appendix}. 

For both Lemma \ref{subpoly counting} and Lemma \ref{quasipoly counting}, we need to choose an appropriate $m = m(n)$ in (\ref{key amalgamation bound}). On one hand, we would like to take $m$ to be sufficiently small for $\binom{m}{2}$ to not be the leading term on the right-hand side of (\ref{key amalgamation bound}), giving $\mathfrak{r}_\mathcal{P}(nm) \leq O\left(\mathfrak{r}_\mathcal{Q}(m)\mathfrak{r}_\mathcal{P}(n)\right)$. On the other hand, $m$ needs to be sufficiently large for us to utilise the asymptotic guarantee on $\mathfrak{r}_\mathcal{Q}(m)$. 

\begin{proof}[Proof of Lemma \ref{quasipoly counting}]
    As $\mathfrak{r}_\mathcal{Q}(n) = e^{O\left(\log^q n\right)}$, there is a constant $C>0$ such that $\mathfrak{r}_\mathcal{Q}(n) \leq e^{C\log^q n}/4 - 1$ for all $n\geq 2$. Let $p = \frac{1}{2 - q}\in [1/2, 1)$ as in the lemma statement. As the amalgamation of a $\mathcal{P}$-colouring and a $\mathcal{Q}$-colouring is always a $\mathcal{P}$-colouring, by Lemma \ref{amalgamation colour count}, we have
    \begin{equation*} 
        \mathfrak{r}_\mathcal{P}(nm) \leq \left(2\mathfrak{r}_\mathcal{Q}(m) + 1\right)\mathfrak{r}_\mathcal{P}(n) + \binom{m}{2}, 
    \end{equation*}
    for any positive integers $m, n$. Therefore, whenever $m\geq 2$, we have
    \begin{equation*} 
        \mathfrak{r}_\mathcal{P}(nm) \leq 1/2\cdot e^{C\log^q m}\mathfrak{r}_\mathcal{P}(n) + \binom{m}{2}. 
    \end{equation*}
    For any positive integer $n\geq 2$, we may take $m = \floor{e^{\log^p n}}\geq 2$,  giving \begin{equation*} 
        \mathfrak{r}_\mathcal{P}(\floor{e^{\log^p n}}n) \leq 1/2\cdot e^{C\log^{pq} n}\mathfrak{r}_\mathcal{P}(n) + 1/2\cdot e^{2\log^p n}. 
    \end{equation*}
    We construct the strictly increasing sequence $n_0, \cdots, n_k, \cdots$ inductively, taking $n_0 = 2$ and $n_{k + 1} = \floor{e^{\log^p n_k}}n_k$ for $k\geq 0$. Therefore, we have
    \begin{equation} \label{quasi-poly lem eq}
        \mathfrak{r}_\mathcal{P}(n_{k + 1}) \leq 1/2\cdot e^{C\log^{pq} n_k}\mathfrak{r}_\mathcal{P}(n_k) + 1/2\cdot e^{2\log^p n_k}. 
    \end{equation}
    Take a constant $\Tilde{C} \geq \max(2, \frac{2C}{2^p - 1})$ such that $\mathfrak{r}_\mathcal{P}(n_0) \leq e^{\Tilde{C}\log^p n_0}$. We prove by induction on $k\geq 0$ that \begin{equation} \label{quasi-poly lem ind eq}
        \mathfrak{r}_\mathcal{P}(n_k) \leq e^{\Tilde{C}\log^{p} n_k}. 
    \end{equation}
    The base case $k = 0$ holds by our choice of $\Tilde{C}$. Inductively, assume that (\ref{quasi-poly lem ind eq}) holds for some $k\geq 0$. By (\ref{quasi-poly lem eq}), we have
    \begin{align*}
        \mathfrak{r}_\mathcal{P}(n_{k + 1}) &\leq 1/2\cdot e^{C\log^{pq} n_k + \Tilde{C}\log^p n_k} + 1/2\cdot e^{2\log^p n_k}\\
        &\leq \max\left(e^{C\log^{pq} n_k + \Tilde{C}\log^p n_k}, e^{2\log^p n_k}\right)\\
        & = e^{C\log^{pq} n_k + \Tilde{C}\log^p n_k}, 
    \end{align*}
    where the last equality follows as $\Tilde{C}\geq 2$. As $n_k\geq n_0 = 2$, we have $\floor{e^{\log^p n_k}} \geq e^{\frac{1}{2}\log^p n_k}$, therefore $n_{k + 1}\geq e^{\frac{1}{2}\log^p n_k}n_k$ and $\log n_{k + 1}\geq \log n_k + \frac{1}{2}\log^p n_k$. Hence, we have
    \begin{equation*}
        \log^p n_{k + 1}\geq \left(\log n_k + \frac{1}{2}\log^p n_k\right)^p \geq \left(1 + \frac{1}{2}\log^{p - 1} n_k\right)^p\log^p n_k. 
    \end{equation*}
    Note that $\frac{1}{2}\log^{p - 1} n_k \leq 1$. For all $x\leq 1$, we have $(1 + x)^p\geq 1 + \left(2^p - 1\right) x$. Therefore, 
    \begin{align*}
        \log^p n_{k + 1}&\geq \left(1 + \frac{1}{2}\log^{p - 1} n_k\right)^p\log^p n_k\\
        &\geq \log^p n_k + \frac{2^p - 1}{2}\cdot \log^{2p - 1} n_k \\
        &= \log^p n_k + \frac{2^p - 1}{2}\cdot \log^{pq} n_k, 
    \end{align*}
    as $2p - 1 = pq$ with $p = \frac{1}{2 - q}$. Therefore, we have 
    \begin{equation*}
        \mathfrak{r}_\mathcal{P}(n_{k + 1}) \leq e^{\Tilde{C}\log^{p} n_{k + 1}}, 
    \end{equation*}
    since $\Tilde{C}\geq \frac{2C}{2^p - 1}$, which completes our inductive proof of (\ref{quasi-poly lem ind eq}). 
    Now for any positive integer $n\geq 3$, take the least $k\geq 0$ such that $n_k\geq n$. Such a $k$ exists as $n_0, \cdots, n_k, \cdots$ is a strictly increasing sequence of positive integers. As $n\geq 3$, we have $k\geq 1$ and therefore $n_{k - 1} < n$. Since $n_k\leq n_{k - 1}^2$, we have $n_k < n^2$. As $\mathfrak{r}_\mathcal{P}$ is an increasing function, we have
    \begin{equation*}
        \mathfrak{r}_\mathcal{P}(n)\leq \mathfrak{r}_\mathcal{P}(n_k)\leq e^{\Tilde{C}\log^{p} n_k}\leq e^{\Tilde{C}\cdot 2^p\log^{p} n}, 
    \end{equation*}
    for all $n\geq 3$. Hence $\mathfrak{r}_\mathcal{P}(n) = e^{O\left(\log^{p} n\right)}$. 
\end{proof}

To apply Lemma \ref{quasipoly counting}, we need to identify appropriate properties $\mathcal{P}$ and $\mathcal{Q}$ of edge-colourings along with suitable values of $p, q$. We shall make the following choices. 
\begin{itemize}
    \item $K_4$-uniqueness (Proposition \ref{k4 prop}). We take $\mathcal{P} = K_4$-uniqueness. We do not need to impose any property $\mathcal{Q}$ upon the auxiliary colouring used. In fact, we may take the auxiliary colouring to be the \emph{trivial colouring}, which assigns the same colour to each edge. We may therefore take $q = 0$ here, resulting in $p = \frac{1}{2}$, as needed. 
    \item $K_5$-uniqueness (Proposition \ref{k5 prop}). We strengthen the inductive hypothesis and define a $\mathcal{P}$-colouring to be an edge-colouring which is both $K_3$-unique (admitting no monochromatic triangles) and $K_5$-unique. Again, we do not need to impose any property $\mathcal{Q}$ upon the auxiliary colouring and the trivial colouring may be used. Once more, we take $q = 0$ and $p = \frac{1}{2}$ here. 
    \item $K_8$-oddness (Theorem \ref{main thm}). We again strengthen the inductive hypothesis and define a $\mathcal{P}$-colouring to be an edge-colouring which is both $K_4$-unique and $K_8$-odd. We impose $\mathcal{Q} = K_4$-uniqueness upon the auxiliary colouring that we use. By Proposition \ref{k4 prop}, we may take $q = \frac{1}{2}$, resulting in $p = \frac{2}{3}$, as needed. 
\end{itemize}

The application of Lemma \ref{quasipoly counting} reduces Proposition \ref{k4 prop}, Proposition \ref{k5 prop} and Theorem \ref{main thm} to the following results on amalgamations, respectively. 

\begin{prop} \label{reduced k4 prop}
    The amalgamation $c\otimes d$ of a $K_4$-unique edge-colouring $c$ on $K_n$ and any edge-colouring $d$ on $K_m$ is a $K_4$-unique edge-colouring of $K_{nm}$. 
\end{prop}

\begin{prop} \label{reduced k5 prop}
    The amalgamation $c\otimes d$ of a $K_3$-unique and $K_5$-unique edge-colouring $c$ on $K_n$ and any edge-colouring $d$ on $K_m$ is a $K_3$-unique and $K_5$-unique edge-colouring of $K_{nm}$. 
\end{prop}

\begin{prop} \label{reduced k8 prop}
    The amalgamation $c\otimes d$ of a $K_4$-unique and $K_8$-odd edge-colouring $c$ on $K_n$ and a $K_4$-unique edge-colouring $d$ on $K_m$ is a $K_4$-unique and $K_8$-odd edge-colouring of $K_{nm}$. 
\end{prop}

The proofs of these propositions will be given in Sections \ref{k4}, \ref{k5}, \ref{k8}, respectively. 

\section{General Observations} \label{elem}

Given an edge-colouring $c: E(K_n)\rightarrow \mathfrak{p}(c)$ of $K_n$, a \emph{weakening} $\Tilde{c}$ of $c$ is an edge-colouring of $K_n$ obtained by post-composing $c$ with a function $f$, namely $\Tilde{c} = f\circ c$. Conversely, we call $c$ a \emph{strengthening} of $\Tilde{c}$. When $f$ is injective, the weakening $\Tilde{c} = f\circ c$ gives an equivalent colouring to $c$. Otherwise, the effect of post-composing by $f$ amounts to identifying colours $k_1, k_2\in \mathfrak{p}(c)$ whenever $f(k_1) = f(k_2)$. 

We make the following simple but important observation regarding the monotonicity of even-chromaticity and unique-chromaticity with respect to the weakening of an edge-colouring. 

\begin{obs} \label{key obs}
    Let $c, \Tilde{c}$ be edge-colourings of $K_n$, where $\Tilde{c}$ is a weakening of $c$. A copy of a graph $H$ in $K_n$ that is unique-chromatic with respect to $\Tilde{c}$ must also be unique-chromatic with respect to $c$. Likewise, a copy of a graph $H$ in $K_n$ that is even-chromatic with respect to $c$ must also be even-chromatic with respect to $\Tilde{c}$. 
\end{obs}

The key example of weakenings we consider is that of the amalgamation $c\otimes d$ into its components. 

\begin{obs} \label{weakening example}
    Let $c\otimes d$ be the amalgamation of edge-colourings $c$ and $d$ (as defined in Section \ref{general}). Then the components $(c\otimes d)_1, (c\otimes d)_2, (c\otimes d)_3, (c\otimes d)_4$ of $c\otimes d$ are weakenings of $c\otimes d$. 
\end{obs}

Consider an amalgamation $c\otimes d$ of a $K_t$-unique (\emph{resp.}~$K_t$-odd) edge-colouring $c$ and any auxiliary edge-colouring $d$. Let $S$ be a copy of $K_t$ in $K_{nm}$ which is not unique-chromatic (\emph{resp.}~is even-chromatic) with respect to $c\otimes d$. The memory of the edge-colouring $c$ via $(c\otimes d)_1$, combined with the special colourings of `vertical' edges via $(c\otimes d)_4$, forces the existence of a `rectangular' structure within $S$ when viewed as a subset of the grid $[n]\times [m] = V(K_{nm})$.

\begin{lem} \label{rect lem}
    Let $c$ be a $K_t$-unique (\emph{resp.}~$K_t$-odd) edge-colouring of $K_n$ and let $d$ be any edge-colouring of $K_m$. Let $S$ be a copy of $K_t$ in $K_{nm}$ (on vertex set $[n]\times [m]$) which is not unique-chromatic (\emph{resp.}~is even-chromatic) with respect to the amalgamation $c\otimes d$. We have
    \begin{enumerate} [label=\alph*)]
        \item $S$ must contain two vertices $(x, y_1), (x, y_2)$ sharing the first coordinate. 
        \item For any two vertices $(x, y_1), (x, y_2)$ in $S$ sharing the first coordinate, $S$ must contain the vertices $(x', y_1), (x', y_2)$ for some $x'\neq x$. 
        \item In particular, $S$ must contain four vertices $(x, y_1), (x, y_2), (x', y_1), (x', y_2)$ forming a `rectangle' in the $[n]\times [m]$ grid. 
    \end{enumerate}
\end{lem}

\begin{proof}
    By Observation \ref{key obs} and \ref{weakening example}, $S$ is not unique-chromatic (\emph{resp.}~is even-chromatic) with respect to any of $(c\otimes d)_1, (c\otimes d)_2, (c\otimes d)_3, (c\otimes d)_4$. 
    \begin{enumerate} [label=\alph*)]
        \item Let $S = \{(\Tilde{x}_1, \Tilde{y}_1), \cdots, (\Tilde{x}_t, \Tilde{y}_t)\}$. For the sake of contradiction, let the first coordinates $\Tilde{x}_1, \cdots, \Tilde{x}_t$ be pairwise distinct. By our assumption that $c$ is a $K_t$-unique (\emph{resp.}~$K_t$-odd) edge-colouring of $K_n$, the copy $\{\Tilde{x}_1, \cdots, \Tilde{x}_t\}$ of $K_t$ in $K_n$ is unique-chromatic (\emph{resp.}~is not even-chromatic) with respect to $c$. Hence the copy $S$ of $K_t$ in $K_{nm}$ is unique-chromatic (\emph{resp.}~is not even-chromatic) with respect to $(c\otimes d)_1$, a contradiction. 
        \item We have $(c\otimes d)_4 ((x, y_1), (x, y_2)) = \{y_1, y_2\}$. As $S$ is not unique-chromatic (\emph{resp.}~is even-chromatic) with respect to $(c\otimes d)_4$, there must be another edge in $E(S)$ receiving the colour $\{y_1, y_2\}$ under the edge-colouring $(c\otimes d)_4$. This concludes the proof as all edges in $K_{nm}$ receiving the colour $\{y_1, y_2\}$ under the edge-colouring $(c\otimes d)_4$ are of the form $(x', y_1)(x', y_2)$ for some $x'\in [n]$. 
        \item Immediate given $(a), (b)$. \qedhere
    \end{enumerate}
\end{proof}

The proof of Lemma \ref{rect lem} captures our general strategy. Given a copy $S$ of $K_t$ in $K_{nm}$ which is not unique-chromatic (\emph{resp.}~is even-chromatic), we iteratively build up more structure in $S$ as a subset of the grid $[n]\times [m] = V(K_{nm})$, from which we then derive contradictions. 

\section{A $K_4$-unique edge-colouring with $n^{o(1)}$ colours} \label{k4}

We now give the proof of Proposition \ref{reduced k4 prop}, which we restate here for clarity. 

\textbf{Proposition \ref{reduced k4 prop}} \textit{The amalgamation $c\otimes d$ of a $K_4$-unique edge-colouring $c$ on $K_n$ and any edge-colouring $d$ on $K_m$ is a $K_4$-unique edge-colouring of $K_{nm}$.}

\begin{proof}
    If $c\otimes d$ is not a $K_4$-unique edge-colouring of $K_{nm}$ (on vertex set $[n]\times [m]$), let $S$ be a copy of $K_4$ in $K_{nm}$ which is not unique-chromatic with respect to the edge-colouring $c\otimes d$. By Lemma \ref{rect lem} $(c)$, $S$ must consist of four vertices $(x_1, y_1), (x_1, y_2), (x_2, y_1), (x_2, y_2)$ forming a `rectangle' in the $[n]\times [m]$ grid. Without loss of generality, we assume $x_1 < x_2$ and $y_1 < y_2$. The edge $(x_1, y_1)(x_2, y_2)$ is the unique edge in $E(S)$ to receive the colour `$+$' under the edge-colouring $(c\otimes d)_3$. Therefore $S$ is unique-chromatic with respect to the weakening $(c\otimes d)_3$ of $c\otimes d$. Hence by Observation \ref{key obs}, $S$ is unique-chromatic with respect to $c\otimes d$, a contradiction. 
\end{proof}

\section{A $K_5$-unique edge-colouring with $n^{o(1)}$ colours} \label{k5}

We now give the proof of Proposition \ref{reduced k5 prop}, which we restate here for clarity. 

\textbf{Proposition \ref{reduced k5 prop}} \textit{The amalgamation $c\otimes d$ of a $K_3$-unique and $K_5$-unique edge-colouring $c$ on $K_n$ and any edge-colouring $d$ on $K_m$ is a $K_3$-unique and $K_5$-unique edge-colouring of $K_{nm}$.}

\begin{proof}
    If $c\otimes d$ is not a $K_3$-unique edge-colouring of $K_{nm}$ (on vertex set $[n]\times [m]$), let $S$ be a copy of $K_3$ in $K_{nm}$ which is not unique-chromatic with respect to $c\otimes d$. By Lemma \ref{rect lem} $(c)$, $S$ must contain at least four vertices, a contradiction. 

    If $c\otimes d$ is not a $K_5$-unique edge-colouring of $K_{nm}$ (on vertex set $[n]\times [m]$), let $S$ be a copy of $K_5$ in $K_{nm}$ which is not unique-chromatic with respect to $c\otimes d$. By Lemma \ref{rect lem} $(c)$, $S$ must contain four vertices $(x_1, y_1), (x_1, y_2), (x_2, y_1), (x_2, y_2)$ forming a `rectangle' in the $[n]\times [m]$ grid. Let $R = \{(x_1, y_1), (x_1, y_2), (x_2, y_1), (x_2, y_2)\}\subseteq S$ and let $(x, y)$ be the fifth vertex in $S$. 

    \textbf{Case 1:} $x_1, x_2, x$ are pairwise distinct. 
    
    As $c$ is $K_3$-unique, there is a colour $k$ occupying exactly one edge $e$ in the copy $\{x_1, x_2, x\}$ of $K_3$ in $K_n$ with edge-colouring $c$. Without loss of generality, we assume $x_1 < x_2$ and $y_1 < y_2$. 

    If $e = x_1x_2$, then the edge $(x_1, y_1)(x_2, y_2)$ is the unique edge in $E(S)$ to receive the colour $k$ under $(c\otimes d)_1$ and the colour `$+$' under the edge-colouring $(c\otimes d)_3$, a contradiction. 

    Let $e = xx_i$ where $i = 1$ or $2$. The edges $(x, y)(x_i, y_1)$, $(x, y)(x_i, y_2)$ are the only edges in $E(S)$ to receive the colour $k$ under $(c\otimes d)_1$. As $S$ is not unique-chromatic, we must have $c\otimes d ((x, y), (x_i, y_1)) = c\otimes d ((x, y), (x_i, y_2))$. In particular, $(c\otimes d)_3 ((x, y), (x_i, y_1)) = (c\otimes d)_3 ((x, y), (x_i, y_2))$. Therefore $y_1 < y_2 < y$ or $y < y_1 < y_2$. 
    
    Take the edges $e_+ = (x_1, y_1)(x_2, y_2), e_- = (x_1, y_2)(x_2, y_1)$ between vertices in $R$. As $S$ is not unique-chromatic, there exist edges $e_+'\neq e_+, e_-'\neq e_-$ in $E(S)$ with $c\otimes d(e_+') = c\otimes d(e_+), c\otimes d(e_-') = c\otimes d(e_-)$. As $e_+$ (\emph{resp.}~$e_-$) is the only edge between vertices in $R$ to receive the colour `$+$' (\emph{resp.}~`$-$') under $(c\otimes d)_3$, the edges $e_+', e_-'$ must be incident to the vertex $(x, y)$. Let $e_+' = (x, y)(x_+, y_+), e_-' = (x, y)(x_-, y_-)$, where $x_+, x_- \in \{x_1, x_2\}, y_+, y_- \in \{y_1, y_2\}$. As $c(x_+, x) = (c\otimes d)_1(e_+') = (c\otimes d)_1(e_+) = c(x_1, x_2)\neq k$, we have $x_+ \neq x_i$. Likewise, $x_- \neq x_i$. Hence $x_+ = x_-$ and $\{y_+, y_-\} = \{y_1, y_2\}$. Therefore, as $y_1 < y_2 < y$ or $y < y_1 < y_2$, we have $(c\otimes d)_3(e_+') = (c\otimes d)_3(e_-')$. However, we have $(c\otimes d)_3(e_+') = (c\otimes d)_3(e_+) = +$ and $(c\otimes d)_3(e_-') = (c\otimes d)_3(e_-) = -$, a contradiction. 

    \textbf{Case 2:} $x_1, x_2, x$ are not pairwise distinct. 
    
    As $x_1\neq x_2$, we may assume that $x = x_1$ without loss of generality. Then $y_1, y_2, y$ are pairwise distinct. Therefore $(x, y)$ is the only vertex in $S$ with second coordinate $y$. Applying Lemma \ref{rect lem} $(b)$ to $(x_1, y_1), (x_1, y)\in S$ gives $(x', y_1), (x', y)\in S$ for some $x'\neq x_1$, a contradiction. 
\end{proof}

\section{A $K_8$-odd edge-colouring with $n^{o(1)}$ colours} \label{k8}

We now give the proof of Proposition \ref{reduced k8 prop}, which we restate here for clarity. 

\textbf{Proposition \ref{reduced k8 prop}} \textit{The amalgamation $c\otimes d$ of a $K_4$-unique and $K_8$-odd edge-colouring $c$ on $K_n$ and a $K_4$-unique edge-colouring $d$ on $K_m$ is a $K_4$-unique and $K_8$-odd edge-colouring of $K_{nm}$.}

\begin{proof}
    By Proposition \ref{reduced k4 prop}, the edge-colouring $c\otimes d$ is $K_4$-unique. 
    
    If $c\otimes d$ is not a $K_8$-odd edge-colouring of $K_{nm}$ (on vertex set $[n]\times [m]$), let $S$ be a copy of $K_8$ in $K_{nm}$ which is even-chromatic with respect to $c\otimes d$. By Lemma \ref{rect lem} $(c)$, $S$ must contain four vertices $(x_1, y_1), (x_1, y_2), (x_2, y_1), (x_2, y_2)$ forming a `rectangle' in the grid $[n]\times [m]$. Let $R = \{(x_1, y_1), (x_1, y_2), (x_2, y_1), (x_2, y_2)\}\subseteq S$. 

    For any $(x, y)\in S\backslash R$, we have $(c\otimes d)_1((x, y), (x_1, y_1)) = (c\otimes d)_1((x, y), (x_1, y_2))$ and $(c\otimes d)_1\linebreak((x, y), (x_2, y_1)) = (c\otimes d)_1((x, y), (x_2, y_2))$. Therefore, in the edge-colouring $(c\otimes d)_1$, every colour occupies an even number of edges among the four edges in $E(\{(x, y)\}, R)$. Summing over $(x, y)\in S\backslash R$, every colour in $(c\otimes d)_1$ occupies an even number of edges among the $16$ edges in $E(S\backslash R, R)$. For edges between vertices in $R$, we have $(c\otimes d)_1((x_1, y_1), (x_2, y_1)) = (c\otimes d)_1((x_1, y_1), (x_2, y_2)) = (c\otimes d)_1((x_1, y_2), (x_2, y_1)) = (c\otimes d)_1((x_1, y_2), (x_2, y_2)) = c(x_1, x_2)$ and $(c\otimes d)_1((x_1, y_1), (x_1, y_2)) = (c\otimes d)_1((x_2, y_1), (x_2, y_2)) = *$. Hence, every colour in $(c\otimes d)_1$ occupies an even number of edges among the $6$ edges in $E(R)$. As $(c\otimes d)_1$ is a weakening of $c\otimes d$, by Observation \ref{key obs}, $S$ is even-chromatic with respect to $(c\otimes d)_1$. Hence every colour in $(c\otimes d)_1$ occupies an even number of edges among the $28$ edges in $E(S)$. Removing $E(S\backslash R, R)$ and $E(R)$ from $E(S)$, we have that every colour in $(c\otimes d)_1$ occupies an even number of edges among the $6$ edges in $E(S\backslash R)$. 

    Let $S\backslash R = \{(\Tilde{x}_1, \Tilde{y}_1), \cdots, (\Tilde{x}_4, \Tilde{y}_4)\}$. If $\Tilde{x}_1, \cdots, \Tilde{x}_4$ are pairwise distinct, then $\{\Tilde{x}_1, \cdots, \Tilde{x}_4\}$ is an even-chromatic copy of $K_4$ in $K_n$ with respect to the $K_4$-unique edge-colouring $c$, a contradiction. 

    Therefore, without loss of generality, let $\Tilde{x}_1 = \Tilde{x}_2$. We now have $(c\otimes d)_1((\Tilde{x}_3, \Tilde{y}_3), (\Tilde{x}_1, \Tilde{y}_1)) = (c\otimes d)_1((\Tilde{x}_3, \Tilde{y}_3), (\Tilde{x}_2, \Tilde{y}_2))$ and $(c\otimes d)_1((\Tilde{x}_4, \Tilde{y}_4), (\Tilde{x}_1, \Tilde{y}_1)) = (c\otimes d)_1((\Tilde{x}_4, \Tilde{y}_4), (\Tilde{x}_2, \Tilde{y}_2))$. Hence each colour in $(c\otimes d)_1$ occupies an even number of the $4$ edges in $E(\{(\Tilde{x}_1, \Tilde{y}_1), (\Tilde{x}_2, \Tilde{y}_2)\}, \{(\Tilde{x}_3, \Tilde{y}_3), (\Tilde{x}_4, \Tilde{y}_4)\})$. Therefore, $(c\otimes d)_1((\Tilde{x}_3, \Tilde{y}_3), (\Tilde{x}_4, \Tilde{y}_4)) = (c\otimes d)_1((\Tilde{x}_1, \Tilde{y}_1), (\Tilde{x}_2, \Tilde{y}_2)) = *$, and thus $\Tilde{x}_3 = \Tilde{x}_4$. 

    We have now partitioned $S$ into $4$ pairs of vertices, each with equal first coordinates. These first coordinates are $x_1, x_2, \Tilde{x}_1, \Tilde{x}_3$. 

    \textbf{Case 1:} $x_1, x_2, \Tilde{x}_1, \Tilde{x}_3$ are pairwise distinct. 

    Here the only edges in $E(S)$ which are coloured $*$ under $(c\otimes d)_1$ are the $4$ edges between each pair of vertices of equal first coordinate. Since $(c\otimes d)_4((x_1, y_1), (x_1, y_2)) = (c\otimes d)_4((x_2, y_1), (x_2, y_2)) = \{y_1, y_2\}$, for $S$ to be even-chromatic, we must have $(c\otimes d)_4((\Tilde{x}_1, \Tilde{y}_1), (\Tilde{x}_2, \Tilde{y}_2)) = (c\otimes d)_4((\Tilde{x}_3, \Tilde{y}_3), (\Tilde{x}_4, \Tilde{y}_4))$. Therefore $\{\Tilde{y}_1, \Tilde{y}_2\} = \{\Tilde{y}_3, \Tilde{y}_4\}$ and we assume, without loss of generality, that $\Tilde{y}_1 = \Tilde{y}_3$ and $\Tilde{y}_2 = \Tilde{y}_4$. Hence $S\backslash R = \{(\Tilde{x}_1, \Tilde{y}_1), (\Tilde{x}_1, \Tilde{y}_2), (\Tilde{x}_3, \Tilde{y}_1), (\Tilde{x}_3, \Tilde{y}_2)\}$ also forms a `rectangle' in the $[n]\times [m]$ grid. We relabel $S\backslash R = \{(x_1', y_1'), (x_1', y_2'), (x_2', y_1'), (x_2', y_2')\}$ to restore the symmetry between $S\backslash R$ and $R$. 

    As $c$ is $K_4$-unique, there is a colour $k$ occupying exactly one edge $e$ in the copy $\{x_1, x_2, x_1', x_2'\}$ of $K_4$ in $K_n$ with edge-colouring $c$. 

    \textbf{Sub-case 1.1:} $e = x_1x_2$ (or likewise if $e = x_1'x_2'$). 

    Without loss of generality, we assume $x_1 < x_2$ and $y_1 < y_2$. The edge $(x_1, y_1)(x_2, y_2)$ is the unique edge in $E(S)$ to receive the colour $k$ under $(c\otimes d)_1$ and `$+$' under $(c\otimes d)_3$, a contradiction. 

    \textbf{Sub-case 1.2:} $e = x_1x_1'$ (or likewise if $e = x_2x_1', x_1x_2', x_2x_2'$). 

    The $4$ edges in $E(\{(x_1, y_1), (x_1, y_2)\}, \{(x_1', y_1'), (x_1', y_2')\})$ are the only edges in $E(S)$ to receive the colour $k$ under $(c\otimes d)_1$. We shall first show that $y_1, y_2, y_1', y_2'$ are pairwise distinct. If not, we may take $y_1 = y_1'$ without loss of generality. For $(x_1, y_1)(x_1', y_1')$ to not be the unique edge in $E(S)$ to receive the colour $k$ under $(c\otimes d)_1$ and the colour `$0$' under $(c\otimes d)_3$, we must have $y_2 = y_2'$. As in Sub-case 1.1, we may now assume $x_1 < x_1'$ and $y_1 < y_2$ without loss of generality, and the edge $(x_1, y_1)(x_1', y_2)$ is the unique edge in $E(S)$ to receive the colour $k$ under $(c\otimes d)_1$ and `$+$' under $(c\otimes d)_3$, a contradiction.  
    
    Hence $y_1, y_2, y_1', y_2'$ must be pairwise distinct. As the edge-colouring $d$ of $K_m$ is $K_4$-unique, there is a colour $k'$ in $d$ occupying exactly one edge $e'$ in the copy $\{y_1, y_2, y_1', y_2'\}$ of $K_4$. 

    \textbf{Sub-sub-case 1.2.1:} $e' = y_1y_2$ (or likewise if $e' = y_1'y_2'$). 

    Without loss of generality, we assume that $x_1 < x_2$ and $y_1 < y_2$. The edge $(x_1, y_1)(x_2, y_2)$ is the unique edge in $E(S)$ to receive the colour $k'$ under $(c\otimes d)_2$ and `$+$' under $(c\otimes d)_3$, a contradiction. 

    \textbf{Sub-sub-case 1.2.2:} $e' = y_1y_1'$ (or likewise if $e' = y_2y_1', y_1y_2', y_2y_2'$). 

    The edge $(x_1, y_1)(x_1', y_1')$ is the unique edge in $E(S)$ to receive the colour $k$ under $(c\otimes d)_1$ and the colour $k'$ under $(c\otimes d)_2$, a contradiction.
    
    \textbf{Case 2:} $x_1, x_2, \Tilde{x}_1, \Tilde{x}_3$ are not pairwise distinct. 

    Here we must have at least $4$ vertices in $S$ sharing the same first coordinate, which we take to be $x\in [n]$. Let these $4$ vertices be $(x, w_1), \cdots, (x, w_4)$. For any unordered pair $\{i, j\}$ of distinct indices in $\{1, 2, 3, 4\}$, we may apply Lemma \ref{rect lem} $(b)$ to vertices $(x, w_i), (x, w_j)$, yielding a pair of vertices $(x_{\{i, j\}}, w_i), (x_{\{i, j\}}, w_j)\in S$ for some $x_{\{i, j\}}\neq x$. For each $i \in \{1, 2, 3, 4\}$, we have $(x, w_i)\in S$ and $(x_{\{i, j\}}, w_i)\in S$ for any $j\neq i$. In particular, there are at least $2$ vertices in $S$ with second coordinate $w_i$. As $S$ consists of $8$ vertices, there are exactly $2$ vertices in $S$ with second coordinate $w_i$. Therefore $x_{\{i, j\}} = x_{\{i, j'\}}$ for any indices $j, j'\neq i$. Hence all $x_{\{i, j\}}$ have the same value, which we take to be $x'\neq x$. Therefore $S = \{(x, w_1), \cdots, (x, w_4), (x', w_1), \cdots, (x', w_4)\}$. 

    As $d$ is $K_4$-unique, there is a colour $k'$ occupying exactly one edge $e'$ in the copy $\{w_1, w_2, w_3, w_4\}$ of $K_4$ in $K_m$ with edge-colouring $d$. Without loss of generality, let $e' = w_1w_2$. Again, without loss of generality, we assume that $x < x'$ and $w_1 < w_2$. The edge $(x, w_1)(x', w_2)$ is the unique edge in $E(S)$ to receive the colour $k'$ under $(c\otimes d)_2$ and `$+$' under $(c\otimes d)_3$, a contradiction. 
\end{proof}

\noindent \textbf{Acknowledgements.} The author is grateful to Oliver Janzer for helpful comments and suggestions. The author would also like to thank the anonymous referees for their thoughtful and constructive comments. 

\bibliographystyle{abbrv}
\bibliography{mybib}

\appendix
\section{Bounds on colour usage} \label{colour counting appendix}

Here we prove Lemma \ref{subpoly counting}, which allows us to prove that $\mathfrak{r}_\mathcal{P}(n)$ is sub-polynomial in $n$ given that $\mathfrak{r}_\mathcal{Q}(n)$ is sub-polynomial for the auxiliary property $\mathcal{Q}$. 

\begin{proof}[Proof of Lemma \ref{subpoly counting}]
    It suffices to show that $\mathfrak{r}_\mathcal{P}(n) = O_{\epsilon}\left(n^\epsilon\right)$ for any $\epsilon > 0$. Since $\mathfrak{r}_\mathcal{Q}(n) = n^{o(1)}$, there exists $N=N(\epsilon)$ such that $\mathfrak{r}_\mathcal{Q}(n) < n^\epsilon/3 - 1$ for any $n\geq N$. We assume, without loss of generality, that $\epsilon < 1$ and $N\geq 2$. As the amalgamation of a $\mathcal{P}$-colouring and a $\mathcal{Q}$-colouring is always a $\mathcal{P}$-colouring, by Lemma \ref{amalgamation colour count}, we have
    \begin{equation*} 
        \mathfrak{r}_\mathcal{P}(nm) \leq \left(2\mathfrak{r}_\mathcal{Q}(m) + 1\right)\mathfrak{r}_\mathcal{P}(n) + \binom{m}{2}, 
    \end{equation*}
    for any positive integers $m, n$. Taking $m = N$, we have
    \begin{equation*} 
        \mathfrak{r}_\mathcal{P}(Nn) \leq 2/3\cdot N^\epsilon \mathfrak{r}_\mathcal{P}(n) + N^2, 
    \end{equation*}
    for any positive integer $n$. We may put this into the form
    \begin{align*} 
        \frac{\mathfrak{r}_\mathcal{P}(Nn)}{(Nn)^\epsilon} &\leq 2/3\cdot \frac{\mathfrak{r}_\mathcal{P}(n)}{n^\epsilon} + N^{2 - \epsilon}\\
        &\leq \max\left(\frac{\mathfrak{r}_\mathcal{P}(n)}{n^\epsilon}, 3N^{2 - \epsilon}\right). 
    \end{align*}
    Hence for any non-negative integer $k$ and any positive integer $n$, we have
    \begin{equation*} 
        \frac{\mathfrak{r}_\mathcal{P}(N^kn)}{(N^kn)^\epsilon} \leq \max\left(\frac{\mathfrak{r}_\mathcal{P}(n)}{n^\epsilon}, 3N^{2 - \epsilon}\right). 
    \end{equation*}
    In particular, taking $n = 1$, we have
    \begin{equation*} 
        \frac{\mathfrak{r}_\mathcal{P}(N^k)}{(N^k)^\epsilon} \leq \max\left(\mathfrak{r}_\mathcal{P}(1), 3N^{2 - \epsilon}\right). 
    \end{equation*}
    Now for any positive integer $n$, we may take $k = \lceil\log_N{n}\rceil$. As $n\leq N^k\leq Nn$ and $\mathfrak{r}_\mathcal{P}$ is an increasing function, we have
    \begin{equation*}
        \mathfrak{r}_\mathcal{P}(n)\leq \mathfrak{r}_\mathcal{P}(N^k)\leq \max\left(\mathfrak{r}_\mathcal{P}(1), 3N^{2 - \epsilon}\right)\cdot (N^k)^\epsilon\leq N\cdot\max\left(\mathfrak{r}_\mathcal{P}(1), 3N^{2 - \epsilon}\right)\cdot n^\epsilon = O_{\epsilon}\left(n^\epsilon\right), 
    \end{equation*}
    as needed. 
\end{proof}

\section{A proof of Proposition \ref{non-clique implies non-unique}} \label{non-clique implies non-unique appendix}

\textbf{Notation.} Whilst a copy of a clique $K_t$ in $K_n$ may be determined by its vertex set, for a general graph $H$, the vertex set of the copy is not sufficient to determine a copy of $H$ in $K_n$. We shall model a copy of $H$ in $K_n$ by an injection $f: V(H)\rightarrow V(K_n)$. For an edge $e = v_1v_2\in E(H)$, we denote by $f(e) = f(v_1)f(v_2)\in E(K_n)$ the corresponding edge in the image. For a graph $G$ and a subset $S\subseteq V(G)$ of vertices, we denote by $G - S$ the graph obtained from $G$ by removing all vertices in $S$. For a graph $G$ and a subset $S\subseteq V(G)$ of vertices, we denote by $N_G(S) = \{v\in V(G)|uv\in E(G)\text{ for all }u\in S\}$ the common neighbourhood of the vertices in $S$. 

We begin with the following lemma, which is the analog of Lemma 2.1 in \cite{versteegen2023upper} for unique-chromaticity in place of even-chromaticity. 

\begin{lem} \label{appendix b lem}
    Let $I$ be an independent set of $t$ vertices in a graph $H$, such that $|E(I, V(H)\backslash I)|\geq 2$. If $u_{H - I}(n) = \Omega\left(n^\epsilon\right)$ for some $\epsilon > 0$, then $u_H(n) = \Omega\left(n^\delta\right)$ where $\delta = \frac{\epsilon}{1 + t\epsilon} > 0$. In particular, if $u_{H - I}(n) = n^{\Omega(1)}$, then $u_H(n) = n^{\Omega(1)}$. 
\end{lem}

\begin{proof}
    Let $u_{H - I}(n) > C \cdot n^\epsilon$ for all $n\geq N$, where $C$ is a positive constant. Let $\Tilde{C} = C^{\frac{1}{1 + t\epsilon}}/2 > 0$. We shall prove that $u_H(n) > \Tilde{C}\cdot n^{\delta}$ for all $n\geq \Tilde{N} = \max\left(3, \left(6\Tilde{C}t\right)^{\frac{1}{1 - \delta}}, \left(\left(2\Tilde{C}\right)^tN\right)^{\epsilon/\delta}\right)$. 
    
    Consider any edge-colouring $c$ of $K_n$ using at most $\Tilde{C}\cdot n^{\delta}$ colours. By the pigeonhole principle, there is a colour $k$ occupying at least $\frac{\binom{n}{2}}{\Tilde{C}\cdot n^{\delta}}$ edges under $c$. Take $G$ to be spanning subgraph of $K_n$ consisting of edges with colour $k$ under $c$. $G$ has average degree at least $\frac{n - 1}{\Tilde{C}\cdot n^{\delta}}\geq \frac{n^{1 - \delta}}{2\Tilde{C}} + t$. Let $U\subseteq V(K_n)$ be a uniformly random subset of $t$ vertices. By linearity of expectation, the expected size of $N_G(U)$ is 
    \begin{equation*}
        \frac{1}{\binom{n}{t}}\sum_{v\in V(K_n)}\binom{\deg_G(v)}{t}. 
    \end{equation*}
    We may extend $\binom{\cdot}{t}$ to a convex function 
    \begin{equation*}
        \binom{x}{t} = \begin{cases}
            x(x - 1)\cdots(x - t + 1)/t! & \text{ if } x\geq t, \\
            0 & \text{ otherwise}
        \end{cases}
    \end{equation*}
    on domain $\mathbb{R}$. By Jensen's inequality, we have
    \begin{equation*}
        \frac{1}{\binom{n}{t}}\sum_{v\in V(K_n)}\binom{\deg_G(v)}{t}\geq \frac{n}{\binom{n}{t}}\binom{\frac{n^{1 - \delta}}{2\Tilde{C}} + t}{t}\geq \frac{n}{n^t}\left(\frac{n^{1 - \delta}}{2\Tilde{C}}\right)^t = (2\Tilde{C})^{-t}\cdot n^{1 - t\delta} = (2\Tilde{C})^{-t}\cdot n^{\delta/\epsilon}. 
    \end{equation*}
    Hence there exists a subset $U$ of $t$ vertices with $|N_G(U)|\geq (2\Tilde{C})^{-t}\cdot n^{\delta/\epsilon}$. Let $W = N_G(U)$ and let $m = |N_G(U)|\geq (2\Tilde{C})^{-t}\cdot n^{\delta/\epsilon}\geq N$. The edge-colouring $c$ of $K_n$ induces an edge-colouring of $K_n[W]$ using at most $\Tilde{C}\cdot n^\delta \leq C\cdot \left((2\Tilde{C})^{-t}\cdot n^{\delta/\epsilon}\right)^\epsilon\leq C\cdot m^\epsilon < u_{H - I}(m)$ colours. Hence there is a copy $f_1: V(H - I)\rightarrow W$ of $H - I$ in $K_n[W]$ which is not unique-chromatic. As $|U| = |I| = t$, there exists a bijection $f_2: I\rightarrow U$. Combining $f_1$ and $f_2$ gives a copy $f: V(H)\rightarrow V(K_n)$ of $H$ in $K_n$. We shall prove that $f$ is not unique-chromatic. 

    If $f$ is unique-chromatic, let $k'$ be a colour occupying exactly one edge $f(e')$ of the copy $f$ of $H$ in $K_n$, where $e\in E(H)$. For any edge $e = v_1v_2\in E(I, V(H)\backslash I)$, let $v_1\in I, v_2\in V(H)\backslash I$. Therefore, $f(v_1)\in U, f(v_2)\in W$. As all edges between $U$ and $W$ have colour $k$, the edge $f(e) = f(v_1)f(v_2)$ corresponding to $e$ must have colour $k$ for any $e\in E(I, V(H)\backslash I)$. As $|E(I, V(H)\backslash I)|\geq 2$, the colour $k$ occupies at least two edges in the copy $f$ of $H$. Hence, $k'\neq k$ and $e'\notin E(I, V(H)\backslash I)$. 

    Since $I$ is an independent set in $H$, we must have $e'\in E(V(H)\backslash I)$. In other words, $f(e')$ is an edge in the copy $f_1$ of $H - I$ in $K_n[W]$, which is not unique-chromatic. Therefore, there exists another edge $e''\in E(V(H)\backslash I)\subseteq E(H)$ such that $f(e'') = f(e') = k'$ under the edge-colouring $c$, a contradiction. Hence $f$ is not unique-chromatic and $u_H(n) > \Tilde{C}\cdot n^{\delta}$ for all $n\geq \Tilde{N}$. 
\end{proof}

\begin{remark}
    In the case that $u_{H - I}(n) = \infty$ (for sufficiently large $n$), which takes place precisely when $H - I$ has no edges, the argument above gives a lower bound of $u_H(n) = \Omega\left(n^{1/t}\right)$. 
\end{remark}

We now show that an appropriate independent set $I$ can always be identified. 

\begin{lem} \label{appendix b lem 2}
    For any non-complete graph $H$ without isolated vertices, there exists an independent set $I$ of vertices in $H$, such that
    \begin{enumerate}[label=\alph*)]
        \item we have $|E(I, V(H)\backslash I)|\geq 2$, 
        \item after removing all isolated vertices in $H - I$, either no vertices remain or we obtain a non-complete graph. 
    \end{enumerate}
\end{lem}

\begin{proof}
    Assume otherwise. Let $S\subseteq V(H)$ be the subset of all vertices of degree at least two in $H$. For any $v\in S$, $I = \{v\}$ satisfies condition $(a)$. Therefore $H - \{v\}$ must consist solely of a clique (of at least two vertices) and isolated vertices. That is, we may partition the vertex set of $H - \{v\}$ into two subsets $C_v$, $I_v$ such that $H - \{v\}$ is complete on $C_v$ and all vertices in $I_v$ are isolated in $H - \{v\}$. Vertices in $I_v$ must neighbour precisely $v$ in $H$. In particular, $S, I_v$ are disjoint and $I_v$'s are mutually disjoint for different choices of $v\in S$. As a result, $H$ is complete on $S\backslash\{v\}\subseteq C_v$. 

    \textbf{Case 1:} $S$ is empty. 

    As all vertices in $H$ have degree exactly one, $H$ must be a perfect matching. As $H$ is not complete, it must consist of at least two edges. Taking $I$ to be the independent set consisting of one endpoint from each edge suffices. 

    \textbf{Case 2:} $1\leq |S|\leq 2$. 

    In this case $H$ cannot have any triangles. Indeed, any vertex in a triangle $T\subseteq V(H)$ has degree at least two. Therefore $T\subseteq S$, contradicting the assumption that $|S|\leq 2$. 

    Take $v\in S$. Recall that $H$ is complete on $C_v$ and $|C_v|\geq 2$. As there are no triangles in $H$, we must have $|C_v| = 2$. Let $C_v = \{u, w\}$. 

    As $uw\in E(H)$ and $\{v, u, w\}$ cannot form a triangle in $H$, we must have $vu\notin E(H)$ or $vw\notin E(H)$. Without loss of generality, we assume that $vu\notin E(H)$. Taking $I = \{u, v\}$ suffices, as $E(H - \{u, v\}) = \emptyset$.
    
    \textbf{Case 3:} $|S|\geq 3$. 
    
    For any pair of distinct vertices $v_1, v_2$ in $S$, we may choose a third vertex $v_3\in S$. Since $H$ is complete on $S\backslash\{v_3\}$, we have $v_1v_2\in E(H)$. Hence, $H$ is complete on $S$. 

    Therefore, as $H$ is not complete, $S$ must be a proper subset of $V(H)$. Take any $u\notin S$. As $u$ has degree one in $H$, there must be a vertex $v\in S$ such that $uv\notin E(H)$. Therefore $u$ has degree one in $H - \{v\}$. In particular $u$ is not isolated in $H - \{v\}$, hence $u\in C_v$. As $S - \{v\}\subseteq C_v$ and $H$ is complete on $C_v$, $S - \{v\}$ must be a subset of the neighbourhood of $u$. This leads to a contradiction, as $|S - \{v\}| \geq 2$ but $u$ has degree one in $H$. 
\end{proof}

We may iterate Lemma \ref{appendix b lem} with the independent set $I$ given by Lemma \ref{appendix b lem 2} to prove Proposition \ref{non-clique implies non-unique}. 

\begin{proof}[Proof of Proposition \ref{non-clique implies non-unique}]
    We induct on the number of vertices $|V(H)|$ of $H$, where $H$ is a non-complete graph without isolated vertices. By Lemma \ref{appendix b lem 2}, there exists an independent set $I$ of vertices in $H$, such that
    \begin{enumerate}[label=\alph*)]
        \item we have $|E(I, V(H)\backslash I)|\geq 2$, 
        \item after removing all isolated vertices in $H - I$, either no vertices remain or we obtain a non-complete graph. 
    \end{enumerate}
    Let $\Tilde{H}$ be the graph obtained by removing all isolated vertices in $H - I$. 
    
    If $V(\Tilde{H})$ is empty, by the remark following Lemma \ref{appendix b lem}, we have $u_H(n) = \Omega\left(n^{1/|I|}\right)$. As $H$ has no isolated vertices, $I$ must be a proper subset of $V(H)$. Therefore, $u_H(n) = \Omega\left(n^{\frac{1}{|V(H)| - 1}}\right)$. 

    Otherwise, $\Tilde{H}$ is a non-complete graph without isolated vertices. As $|E(I, V(H)\backslash I)|\geq 2$, $I$ must be non-empty. Therefore, $|V(\Tilde{H})| < |V(H)|$. By the inductive hypothesis, we have $u_{\Tilde{H}}(n) = \Omega\left(n^{\frac{1}{|V(\Tilde{H})| - 1}}\right)$. Since the addition of isolated vertices to a graph $H'$ does not change the function $u_{H'}$, we have $u_{H - I}(n) = \Omega\left(n^{\frac{1}{|V(\Tilde{H})| - 1}}\right) = \Omega\left(n^{\frac{1}{|V(H)| - |I| - 1}}\right)$. Therefore, by Lemma \ref{appendix b lem}, we have $u_H(n) = \Omega\left(n^{\frac{1}{|V(H)| - 1}}\right)$. 
\end{proof}

\end{document}